\documentclass[11pt]{amsart}
\usepackage{amsmath,amsfonts,latexsym,amssymb,amscd, graphicx,pictexwd,color}

\renewcommand{\P}{{\mathbb P}}
\newcommand{\E}{{\mathbb E}}
\newcommand{\RR}{{\mathbb R}}

\newcommand{\NN}{{\mathbb N}}

\newcommand{\bP}{{\mathbf P}}

\newcommand{\Var}{{\rm Var}}
\newcommand{\Cov}{{\rm Cov}}

\def\iD{{\mathit D}}

\newtheorem{lemm}{Lemma}[section]

\numberwithin{equation}{section}

\def\D{\mathcal{D}}
\def\bP{\mathbf{P}}

\newtheorem{theo}{Theorem}
\newtheorem{coro}{Corollary}

\newenvironment{dem}{\vskip 2mm\noindent {\it Proof} :}
                    {\hfill $\square$ \vskip 2mm \noindent}

\newenvironment{demo}{\vskip 2mm\noindent {\bf }}
                    {\hfill $\square$ \vskip 2mm \noindent}

\begin{document}

\title[Shock fluctuations for the Hammersley process]{Shock fluctuations for the Hammersley process}

\author{Leandro P. R. Pimentel and Marcio W. A. de Souza}
\address{}
\email{}
\thanks{M.W.A. de Souza was partially supported by FAPESP - Fundacao de Amparo a Pesquisa do Estado de Sao Paulo}


\begin{abstract}
We consider the Hammersley interacting particle system starting from a shock initial profile with densities $\lambda,\rho\in\RR$ ($\lambda > \rho$). The microscopic shock is taken as the position of a second-class particle initially at the origin, and the main results are: (i) a central limit theorem for the shock; (ii) the variance of the shock equals $2[\lambda\rho(\lambda - \rho)]^{-1}t + O(t^{2/3})$. By using the same method of proof, we also prove similar results for first-class particles.
\end{abstract}

\maketitle

\section{Introduction}


The Hammersley particle system is a Markov process in the space of counting measures $\nu$ on $\RR$. The atoms of $\nu$ are called particles, and each particle jumps to the left, at rate equal to the distance to the left neighbor, with the new location chosen uniformly at random between the jumper and its left neighbor. The unique family of ergodic invariant measures for the Hammersley process is given by Poisson processes of intensity $\lambda>0$. This particle system was introduced by Aldous and Diaconis \cite{AD}, based on an approach developed by Hammersley \cite{Ha} to study the famous Ulam's problem on the longest increasing subsequence of a random permutation.
\newline


Let $M^t_\nu$ denote the particle configuration at time $t\geq 0$, with $M^0_\nu=\nu$. Consider two initial measures $\nu_1$ and $\nu_2$ such that $\nu_1\leq \nu_2$, i.e. $\nu_1(A)\leq \nu_2(A)$ for every $A\subseteq \RR$. We note that every particle (atom) of $\nu_1$ is also a particle in $\nu_2$, and there may be particles that are in $\nu_2$ but not in $\nu_1$ (discrepancies). If one runs simultaneously (basic coupling) $M^t_{\nu_1}$ and $M^t_{\nu_2}$ then $M^t_{\nu_1}\leq M^t_{\nu_2}$ for all $t\geq 0$.  This property gives rise to the so called two-class Hammersley process, which takes into account the evolution of the particles in $\nu_1$, called \emph{first-class particles}, and the evolution of the discrepancies between $\nu_1$ and $\nu_2$, called \emph{second-class particles}. In the present paper we study the fluctuations of first-class and second-class particles, as well as the flux of these particles.      
\newline 
 

The hydrodynamic limit of the density of particles in the Hammersley process is described by the Burgers equation \cite{AD,SeM}: $\partial_t u +\partial_x f( u) =0$, with $f(u)=-1/u$. The intuition behind it relies on the idea that the time derivative of the flux of particles (at position $x$ at time $t$) should be approximately the reciprocal of the space derivative of the flux (see equation (9) in \cite{AD}). Thus, if one sets $u=\partial_x U$, then $U$ solves $\partial_t U + f(\partial_x U)=0$. The solutions of the Burgers equation may develop traveling shocks (discontinuities), and this behavior is due to  collisions between characteristic curves emanating from different positions. A simple example is given by the initial profile $u(x)=\rho$ for $x<0$ and $u(x)=\lambda$ for $x>0$, where $\lambda>\rho$. For this initial profile the characteristic curves are straight lines: $z(t)=x_0+f'(u(x_0))t=x_0+u(x_0)^{-2}t$. The solution is $u(x,t)=u(x-ct)$, where $c$ is the velocity of the shock. To find out the value of $c$, one must notice that (for $a<ct<b$),
$$\int_a^b\partial_t u(x,t)dx= f(u(a,t))-f(u(b,t))= f(u(a-ct))-f(u(b-ct))=-\rho^{-1}+\lambda^{-1}\,,$$ 
and that
$$\int_a^bu(x,t)dx= \int_a^{ct} u(x,t)dx+\int_{ct}^b u(x,t)dx=\rho(ct-a)+\lambda(b-ct)\,,$$
which yields to $c=(\lambda\rho)^{-1}$.
\newline

At the microscopic level, a shock related to $(M^t_\nu)_{t \geq 0}$ is a random position (depending on $t$) with the property that (uniformly in time) the asymptotic densities to the right and left of the shock are different. Typically, a shock can be described by a second-class particle \cite{Fe,Fe2}. Ferrari and Fontes \cite{FeFoS} is the main reference on the fluctuations of a microscopic shock for the exclusion process. In the asymmetric simple exclusion process, they obtained a central limit theorem and the diffusion coefficient for the position of a second class particle initially at the origin (in the shock regime). For the Hammersley process, Cator and Groeneboom \cite{CG1} proved a law of large numbers for the second class particle and Sepp\"al\"ainen \cite{SeD} proved functional central limit theorems for the position of a first-class tagged particle for a large class of shock initial profiles. 
\newline 

The main results of this paper are: (i) a central limit theorem for the location of a second-class particle starting from a shock profile $\nu_{\rho,\lambda}$, composed by two independent Poisson measures with densities $\rho$ and $\lambda$, with $\rho<\lambda$, to the left and to the right of the origin, respectively; (ii) the variance of the second-class particle equals $2[\lambda\rho(\lambda - \rho)]^{-1}t + O(t^{2/3})$. The $t^{2/3}$ term comes from the fluctuations of the dynamics in equilibrium (along the characteristic direction), which are known to be of this order since the seminal work of Cator and Groeneboom \cite{CG2} (see also \cite{BaRa}). Our method of proof is based on a combination of this cube-root asymptotic \cite{CG2}, together with the analysis of the influence of the initial condition on the flux \cite{CPS}. To compute the variance, we use a key result that allows us to relate moments of the flux of particles with moments of the location of a tagged particle (Lemma \ref{lemm:relation} and Lemma \ref{lemm:relation2}). Although our techniques have some similarities with the methodology developed by Ferrari and Fontes \cite{FeFoS}, some extra effort is needed in order to obtain the cube root error (Theorem \ref{Theo:2fluxVar}). Furthermore, we have included fluctuations results in the stationary regime of the two-class model, which were not considered so far. In section 2 we define the model and state the results, and in section 3 we prove them.
\newline

\section{Results: central limit theorems and diffusion coefficients}

\subsection{First-class particles} 
We can formally define the Hammersley process by giving its graphical construction. Let $\nu$ be a counting measure on $\RR$ with a positive density of particles. We write $z\in\nu$ whenever $\nu(\{z\})=1$. One can think of it as the initial profile of particles and the dynamics have the following rule: let $\bP\subseteq\RR^2$ be a two-dimensional homogeneous Poisson random set of intensity one, that is assumed to be independent of $\nu$, whenever $\nu$ is random. Set $M_\nu^0=\nu$ and if there is a Poisson point at $(x,t)$ then define 
$$M^t_\nu(\{x\})=M^{t^{-}}_\nu(\{x\})+1\,,$$
and for $y>x$ 
\begin{equation}\nonumber\label{eq:dynamics}
M^t_\nu((x,y])=\left(M^{t^{-}}_\nu((x,y])-1\right)_+\,.
\end{equation}
Here, $M^{t^{-}}_\nu$ is the configuration of particles at time $t$ if the Poisson point at $(x,t)$ would be removed. To the left of $x$ the measure does not change. In words, for $(x,t)\in\bP$, at time $t> 0$ the nearest particle to the right of $x$ is moved to $x$. Denote by $\nu_\lambda$ a Poisson (process) counting measure on $\RR$ of intensity $\lambda>0$. Then $\nu_\lambda$ is invariant for the evolution \cite{AD}:
$$M^t_{\nu_\lambda}\stackrel{\D}{=}\nu_\lambda\,\,\mbox{ for all $t\geq0$}\,.$$ 

Let $L_\nu(x,t)$ denote the flux of particles through the space-time line connecting $(0,0)$ to $(x,t)$. That is, if we label particles according to their initial positions, $X^z_\nu(0)=z$ for $z\in\nu$, then 
\begin{equation*}
L_\nu(x,t):=\#\left\{ z\in\nu\,:\,z>0\,,\,X^z_\nu(t)\leq x\right\}-\#\left\{ z\in\nu\,:\,z\leq 0\,,\,X^z_\nu(t)> x\right\}\,,
\end{equation*}
and 
\begin{equation*}\nonumber
M^t_\nu((x,y])=L_\nu(y,t)-L_\nu(x,t)\,.
\end{equation*}
Notice that, $L_\nu(x,0)=\nu(x)$ where
\begin{equation}\label{Eq:nu} 
\nu(x):=\,\left\{\begin{array}{ll} \nu((0,x]) & \mbox{for } x> 0\,,\\
0 & \mbox{for } x = 0\,, \\
-\nu((x,0]) & \mbox{for }x <0\,.\end{array}\right.
\end{equation}
Denote by $X_\nu(t)$ the position at time $t$ of the particle that was initially located at the first point of $\nu$ to the left of the origin, that is, $X_\nu(0) := \sup\{z \in \nu : z \leq 0 \}$. The following relation is easily derived from the definition of the flux: for all $t\geq 0$ and $x\in\RR$ we have that
\begin{equation}\label{Eq:FluxL}
 L_{\nu}(x,t) = \,\left\{\begin{array}{ll}  - M_{\nu}^{t}( (x, X_\nu(t)] ) & \mbox{if } x < X_\nu(t)\,, \\
 0 & \mbox{if } x = X_\nu(t)\,, \\
 M_{\nu}^{t}( (X_\nu(t), x] ) & \mbox{if } x > X_\nu(t)\,. \end{array}\right.
\end{equation}
The flux has also a variational representation \cite{AD} given by
\begin{equation}\label{last-evo}
L_{\nu}(x,t):= \sup_{z\leq x} \left\{ \nu(z) + L((z,0),(x,t))\right\}\ \ \ \ \ \mbox{ for }\,\,x\in\RR\,,\mbox{ and }\,\,t\geq 0\,,\\ 
\end{equation}
where $L((0,z),(x,t))$ is the so called last-passage percolation time between $(z,0)$ and $(x,t)$, defined as the maximal length among all increasing sequences (in the two-dimensional sense) of Poisson epochs $\bP$ lying in the rectangle $(z,x]\times(0,t]$.
\newline

Add a particle to $\nu$ at position $x=0$, $\nu'=\nu+\delta_0$, and take $\nu=\nu_\lambda$. The \emph{tagged first-class particle} process $(X_{\nu'_\lambda}(t),t\geq 0)$ is the process that keeps track of the position at time $t$ of the particle which was initially at $0$.

\begin{theo}\label{Theo:distX} 
Let $(P_1(x)\,:\,x\geq 0)$ and $(P_2(t)\,:\,t\geq 0)$ be independent Poisson processes of intensity $\lambda$ and $1/\lambda$, respectively. Then
\begin{equation}\label{1dist}
\P\left(X_{\nu'_\lambda}(t)> x \right) = \P\left(P_1(-x)+1>P_2(t)\right)\,,\mbox{ for $t\geq 0$ and $x< 0$ }\,.
\end{equation}
\end{theo}

\begin{coro}
Let $\Phi$ denote the standard normal distribution function. Define
$$\mu_1:=-\frac{1}{\lambda^2}\,\,\mbox{ and }\,\,\sigma_1:=\sqrt{\frac{2}{\lambda^3}}\,.$$
Then, for all $u \in\RR$,
\begin{equation}\label{1clt}
\lim_{t\to\infty}\P\left(X_{\nu'_\lambda}(t)\leq \mu_1 t + (\sigma_1\sqrt{t})u \right) =\Phi(u)\,.
\end{equation}
\end{coro}
\begin{dem}
Let $z_u(t):=\mu_1 t +(\sigma_1 \sqrt{t}) u$. For fixed $u\in\RR$, we have that $z_u(t)<0$ for all large enough $t$. Hence, by \eqref{1dist}, it is enough to show that  
\begin{equation*}
\lim_{t\to\infty}\P\left(P_1(-z_u(t)) +1 \leq P_2(t)\right) =\Phi(u)\,,
\end{equation*}
which is a straightforward consequence of the central limit theorem for Poisson processes.

\end{dem}

In principle, one should be able to compute the mean and the variance of $X_{\nu'_\lambda}$ with \eqref{1dist} in hands. However we will develop a different argument for these computations, whose ideas will also be useful to deal with second-class particles.    

\begin{theo}\label{Theo:cltX} 
For all $t \geq 0$,
\begin{equation}\label{1mean&var}
\E X_{\nu'_\lambda}(t) = \mu_1 t\,\,\,\mbox{ and }\,\,\,\Var X_{\nu'_\lambda}(t) = \sigma_1^2 t\,.
\end{equation}
\end{theo}

\bigskip
\subsection{Second-class particles}
The Hammersley particle system is attractive: if $\nu((x,y])\leq \bar\nu((x,y])$ for all $x<y$ in $\RR$, and $(M^t_{\nu},M^t_{\bar\nu})$ evolves according to the same Poisson epochs $\bP$ (basic coupling), then 
$$M^t_{\nu}((x,y])\leq M^{t}_{\bar\nu}((x,y])\mbox{ for all $x<y$ in $\RR$}\,.$$ 
Recall that $\nu'$ is obtained from $\nu$ by adding a particle at position $x=0$ and running simultaneously the Hammersley processes $M^t_\nu$ and $M^t_{\nu'}$. In this way, one can keep track of the location of the discrepancy initially at $x=0$ through 
$$Z_{\nu'}(t):=\inf\left\{ x\geq 0\,:\,M^t_{\nu'}=M^t_{\nu}+1\right\}\,.$$
We note that $Z_{\nu'}(t)$ is a non-decreasing function of $t$, meaning that it moves to the right. This discrepancy successively jumps to the previous positions of Hammersley particles directly to the right of it, at times where these particles jump to a position to the left of the discrepancy. This can be seen as a priority rule that subordinates the dynamics of the discrepancy initially at $0$, called the \emph{second-class particle}, to the dynamics of the Hammersley first-class particles.
\newline

Given a Poisson counting measure $\nu_\lambda$ on $\RR$, and for a fixed $\rho<\lambda$, delete the particles of $\nu_\lambda$ with probability $1-\rho\lambda^{-1}$, to obtain a Poisson counting measure $\nu_\rho\leq \nu_\lambda$ on $\RR$ of intensity $\rho>0$. Thus, we can write 
$$\nu_\lambda=\nu_\rho+\bar\nu_{\lambda-\rho}\,,$$
where $\bar\nu_{\lambda-\rho}$ is a Poisson process of intensity $\lambda-\rho$ and independent of $\nu_\rho$. In this way, the process composed by the discrepancies, i.e. 
$$\bar M^t_{\bar\nu_{\lambda-\rho}}:= M^t_{\nu_\lambda}-M^t_{\nu_\rho}\,,$$
can be seen as a particle process as well. We call first-class particles the atoms of $M^t_{\nu_\rho}$ and we call second-class particles the atoms of $\bar{M}^t_{\bar\nu_{\lambda-\rho}}$. 
\newline

Add a discrepancy at $0$ to $\bar\nu_{\lambda-\rho}$, and run both process $M^t_{\nu'_\lambda}$ and $M^t_{\nu_\rho}$ with the same epochs.  We note that the motion of a second class particle is not affected by second class particles which are located at its left, and it is equally affected by first and second class particles to its right \cite{Fe2,FeMaH}. Thus, considering both first and second class particles on its right side and only first class particles on its left side we obtain a Hammersley process with shock initial measure $(M^t_{\nu'_{\rho,\lambda}})_{t \geq 0}$ from a two-class process like $(M^t_{\nu_\rho}, \bar{M}^t_{\bar\nu_{\lambda-\rho}})_{t \geq 0}$. In this context, the \emph{tagged second class particle} which was initially at $0$ is given by $Z_{\nu'_{\rho,\lambda}}(t)$, where $\nu'_{\rho,\lambda}:=\nu_{\rho,\lambda}+\delta_0$ and 
\begin{equation}\label{shock1}
\nu_{\rho,\lambda}(y):=\,\left\{\begin{array}{ll} \nu_\lambda((0,y]) & \mbox{for } y> 0\,,\\
0 & \mbox{for } y = 0\,, \\
-\nu_\rho((y,0]) & \mbox{for }y <0\,.\end{array}\right.
\end{equation}
Let $\pi_x \nu$ be the measure $\nu$ translated by $x \in \RR$. We call by \emph{the process $(M^t_{\nu_\rho}, \bar{M}^t_{\bar\nu_{\lambda-\rho}})_{t \geq 0}$ as seen from the second-class particle $Z_{\nu'_{\rho,\lambda}}(t)$} as the process $(\pi_{Z_{\nu'_{\rho,\lambda}}(t)}M^t_{\nu'_{\rho,\lambda}})_{t \geq 0}$ defined as the translation by $Z_{\nu'_{\rho,\lambda}}(t)$ of $(M^t_{\nu'_{\rho,\lambda}})_{t \geq 0}$, where the last  process is the Hammersley process with shock initial profile induced by the two-class process and $Z_{\nu'_{\rho,\lambda}}(t)$. This construction of the Hamersley process with shock initial measure (from a two class system) will be crucial to analyze the asymptotic behavior of the second-class particle.   


\begin{theo}\label{Theo:CLTShock}
Let 
$$\mu_2:=\frac{1}{\lambda\rho}\,\,\,\,\mbox{ and }\,\,\,\,\sigma_2:=\sqrt{\frac{2}{\lambda\rho(\lambda-\rho)}}\,\,\,.$$
Then, for all $u \in\RR$,
\begin{equation}\label{2clt}
\lim_{t\to\infty}\P\left(Z_{\nu'_{\rho,\lambda}}(t)\leq \mu_2 t +(\sigma_2\sqrt{t}) u\right)=\Phi(u)\,.
\end{equation}
\end{theo}

\begin{theo}\label{Theo:Shock}
For all $t \geq 0$,
\begin{equation}\label{2mean&var}
\E Z_{\nu'_{\rho,\lambda}}(t) = \mu_2 t \,\,\,\mbox{ and }\,\,\,\Var Z_{\nu'_{\rho,\lambda}}(t) = \sigma_2^2 t + O(t^{2/3})\,.
\end{equation}
\end{theo}

\bigskip

\subsection{Flux of second-class particles}
Recall \eqref{last-evo} and define the exit point processes 
$$Y^s_{\nu}(x,t):=\sup\left\{z\leq x\,:\, \nu(z) + L((z,0),(x,t))=L_\nu(x,t)\right\}\,,$$
and 
$$Y^i_{\nu}(x,t):=\inf\left\{z\leq x\,:\, \nu(z) + L((z,0),(x,t))=L_\nu(x,t)\right\}\,.$$
Since the seminal paper \cite{CG2}, they became an important tool to analyze fluctuations of the flux.  

In Corollary 1.2 \cite{CPS} it was proved that the flux of (first-class) particles can be approximated, on the $\sqrt{t}$ scale, by a function of the initial profile. Precisely, for any $x\in\RR$ and $t\geq 0$
\begin{equation}\label{1classInf}
\E\left(\left\{L_{\nu_\lambda}(x,t)-\Big[\nu_\lambda(x-\lambda^{-2}t)+2\lambda^{-1}t\Big]\right\}^2\right)=\Var L_{\nu_\lambda}(\lambda^{-2}t,t)\,.
\end{equation}
Cator and Groeneboom \cite{CG2} prove that 
\begin{equation}\label{cuberoot}
\Var L_{\nu_\lambda}(\lambda^{-2}t,t)=\lambda\E Y^s_{\nu_\lambda}(\lambda^{-2}t,t)_+\sim O(t^{2/3})\,,
\end{equation}
and hence,
$$\E\left(t^{-1}\left\{L_{\nu_\lambda}(x,t)-\Big[\nu_\lambda(x-\lambda^{-2}t)+2\lambda^{-1}t\Big]\right\}^2\right)\sim O(t^{-1/3})\,.$$\\

For the two-class process $(M_{\lambda},\bar{M}_{\rho,\lambda})$ previously introduced, we label particles according to their initial position, $Z^z_{\bar\nu_{\lambda - \rho}}(0)=z$ for $z\in \bar\nu_{\lambda - \rho}$ and define the flux of second-class particles as   
\begin{equation*}
\bar L_{\rho,\lambda}(x,t) :=\#\left\{ y\in \bar\nu_{\lambda - \rho}\,:\,y>0\,,\,Z^y_{\bar\nu_{\lambda - \rho}}(t)\leq x\right\}-\#\left\{ y\in\bar\nu_{\lambda - \rho}\,:\,y\leq 0\,,\,Z^y_{\bar\nu_{\lambda - \rho}}(t)> x\right\}\,,
\end{equation*}
By mass conservation, we have
$$\bar L_{\rho,\lambda}(x,t) = L_{\nu_\lambda}(x,t)- L_{\nu_\rho}(x,t)\,.$$

In the same way to prove \eqref{1classInf}, and also using \eqref{cuberoot}, we get the influence of the initial profile on the flux of second-class particles. Let 
$$\Gamma_{\rho,\lambda}(x,t):=\nu_\lambda(x-\lambda^{-2}t)-\nu_\rho(x-\rho^{-2}t)+2(\lambda^{-1}-\rho^{-1})t\,.$$
By translation invariance of $\nu_\lambda$ and $\nu_\rho$ (see Theorem 1 \cite{CPS}), 
$$\left( L_{\nu_\lambda}(x,t)-\nu_\lambda(x-Vt)\,,\,L_{\nu_\rho}(x,t)-\nu_\rho(x-Vt)\right)\,{\buildrel \iD \over =}\,\left(L_{\nu_\lambda}(Vt,t)\,,\,L_{\nu_\rho}(Vt,t)\right)\,,$$
and thus (take $V=\rho^{-2}$),
$$ \bar L_{\rho,\lambda}(x,t)-\big[ \nu_\lambda(x-\lambda^{-2}t)-\nu_\rho(x-\rho^{-2}t)\big]\,{\buildrel \iD \over =}\,\big[L_{\nu_\lambda}(\rho^{-2}t,t)-\nu_\lambda\left((\rho^{-2}-\lambda^{-2})t\right)\big]-L_{\nu_\rho}(\rho^{-2}t,t)\,.$$
Together with \eqref{1classInf} and \eqref{cuberoot}, this implies that 
\begin{eqnarray}
\nonumber R_1(t) &:=& \E\left(\left\{\bar L_{\rho,\lambda}(x,t)-\Gamma_{\rho,\lambda}(x,t)\right\}^2\right)\\
\nonumber&\leq& 2\Big(\Var L_{\nu_\lambda}(\lambda^{-2}t,t)+\Var L_{\nu_\rho}(\rho^{-2}t,t)\Big)\\
\label{Theo:Flux2}&\sim& O(t^{2/3})\,.
\end{eqnarray}
We emphasize that \eqref{Theo:Flux2} is function of $t$ only. This will be crucial to estimate the variance $\bar L_{\rho,\lambda}(x,t)$ for large values of $x$. 
\begin{theo}\label{Theo:2fluxVar}
For $x\geq \rho^{-2}t$, we have that
\begin{equation}\label{2fluxVar}
\Var\bar L_{\rho,\lambda}(x,t)=(\lambda-\rho)\left(x+\frac{t}{\lambda\rho}\right)+R_1(t)+2\Big(R_2(t)+R_3(x,t)\Big)\,,
\end{equation}
where 
$$R_2(t):=-\rho\E\left(Y^s_{\nu_\rho}(\rho^{-2}t,t)_+\wedge(\rho^{-2}t-
\lambda^{-2}t)\right)\,,$$
and
$$R_3(x,t):=-(\lambda-\rho)\E\Big(Y^i_{\nu_\lambda}(\lambda^{-2}t,t)_-\wedge \left(x-\lambda^{-2}t\right)\Big)-\rho\E\Big( Y^i_{\nu_\lambda}(\lambda^{-2}t,t)_-\wedge \left[(\rho^{-2}-\lambda^{-2})t\right] \Big)\,.$$
\end{theo}

As a consequence of \eqref{cuberoot}, \eqref{Theo:Flux2} and \eqref{2fluxVar}, we have 
\begin{equation}\label{2fluxAsymp}
\exists\,\,\lim_{x\to\infty}\left\{\Var\bar L_{\rho,\lambda}(x,t)-(\lambda-\rho)\left(x+\frac{t}{\lambda\rho}\right)\right\}=R(t)\sim O(t^{2/3})\,.
\end{equation}

\bigskip

\subsection{Stationary multi-class process}\label{twoclass}
 Besides the measure $\nu_{\rho,\lambda}$ \eqref{shock1} be taken as the usual example of shock measure, it is important to stress that this measure is not invariant from the point of view of a second class particle $Z_{\nu'_{\rho,\lambda}}$ added to the origin, that is, the process $(\pi_{Z_{\nu'_{\rho,\lambda}}(t)}M^t_{\nu'_{\rho,\lambda}})_{t \geq 0}$ is not stationary (see the discussion in the end of page $4$). Little was known until the work of Ferrari and Martin \cite{FeMaH}, where they showed that the invariant measure of multi-class versions of the exclusion and Hammersley processes can be seen as the output process of a M/M/1 queue. We denote this measure as $(\nu_\rho, \bar{\nu_\xi})$. Using this relation, they prove that starting a two-class Hammersley process with the measure $(\nu_\rho, \bar{\nu_\xi})$ conditioned to have a second class particle $Z_{\xi'} (t)$ at the origin, then the process as seen from $Z_{\xi'} (t)$ is stationary. 
\newline

 In the previous section, the starting measure was taken $(\nu_\rho,\bar\nu_{\lambda-\rho})$, where $\nu_\rho$, and $\bar\nu_{\lambda-\rho}$ were independent Poisson processes of intensities $\rho$ and $\lambda-\rho$, respectively. We prove that the second-class particle is indeed a microscopic shock for both initial measures, and that its law is affected only by an error of order $1$, when changing the initial measure from $(\nu_\rho,\bar\nu_{\lambda-\rho})$ to $(\nu_\rho,\bar\nu_\xi)$. This implies the same asymptotic variance for the shock when the coupled process is in equilibrium.
\newline

We start describing the Hammersley process with two classes of particles $(M^t_{\nu_{\rho}},\bar M^t_\xi)_{t \geq 0}$ in its stationary regime, using the Ferrari-Martin construction. The invariant measure of the Hammersley process with two classes can be seen as the output process of a stationary $M/M/1$ queue. The stationary $M/M/1$ queue is a well known Markov process that can be constructed as a function of two independent Poisson processes, the arrival process $(A(s))_{s\in\RR}$ of rate $\mu_1$ and the service process $(S(s))_{s\in\RR}$ of rate $\mu_2$, with $\mu_1 < \mu_2$. In fact, it is possible to construct the whole process using the queue function in a coupling called the multi-line process, which is a process composed of $n$ copies of Hammersley processes (different from the basic coupling). For simplicity, we fix $n=2$.  
\newline

Let $\alpha_1$ and $\alpha_2$ denote two homogeneous one-dimensional Poisson processes. We couple two Hammersley processes $(M^t_{\alpha_1},M^t_{\alpha_2})$ as follows: the jumps of the second marginal are determined by the Poisson process $\bP$ as usual; meanwhile, the jumps of the first marginal are determined by the dual process $\bP^*$. The dual points of $\bP^*$ are the positions where there was a particle of the process $(M^t_{\alpha_2})_{t \geq 0}$ immediately before this particle jumps to a point of $\bP$. So, $\bP^*=\bP^*(\bP,\alpha_2)$. Cator and Groeneboom \cite{CG1} showed a Burke's theorem for the Hammersley process, implying  that $\bP^*$ is also a two dimensional Poisson process with the same intensity as $\bP$. In particular, this implies $(M^t_{\alpha_1})_{t\geq 0}$ will also be a Hammersley process. This coupling is called Multi-line (two-line) coupling. Observe that here the initial condition $(\alpha_1, \alpha_2)$ could be taken as any pair of point processes, with no assumptions such as $\alpha_1 \leq \alpha_2$ as we did in the basic coupling before. It can be shown that if $(\alpha_1,\alpha_2)$ are the product of Poisson processes, then the two-line process is stationary and this is the unique ergodic invariant measure \cite{FeMaH}. 
\newline 

To define the queueing system, we take $\alpha_1$ and $\alpha_2$ independent Poisson processes of rates $\rho$ and $\lambda$, with $\rho < \lambda$, respectively. For $t\geq 0$, we construct the queue process $(Q^{t}(x))_{x\in\RR}$ from the processes $M^t_{\alpha_1}$ and $M^t_{\alpha_2}$ by setting $M^t_{\alpha_1}$ as the arrival process and $M^t_{\alpha_2}$ as the service process (observe that time for the queue system is space for the Hammersley process). So, $Q^t$ is a (deterministic) function of the two-line process. Burke's theorem for $M/M/1$ queue says that the departures of the queue (the effective services) form a Poisson process in $\RR$ of rate $\rho$, which we denote $D^t_\rho$. Denote the unused services of the queue by $U^t$ and note that, by definition, $D^t_\rho \subseteq M^t_{\alpha_2}$, $U^t\subseteq M^t_{\alpha_2}$ and $M^t_{\alpha_2}=D^t_\rho\cup U^t$. Thus, $(D^t_\rho,U^t)$ is a function of $Q^t$, and therefore also a function of $(M^t_{\alpha_1},M^t_{\alpha_2})$. The key observation in \cite{FeMaH} was that $(D_\rho^t,U^t)$ defines a  two-class process started from $\nu_{\rho} = D_\rho^0$ and $\bar\nu_\xi = U^0$ and its distribution is invariant for the dynamics. (In this setting, $\nu_\lambda=\alpha_2$.) 
\newline

A process $(Z(t))_{t\geq 0}$ is a microscopic shock for the particle system $(M^t_\nu)_{t \geq 0}$ if, uniformly in t, the asymptotic densities to the right and left of $Z(t)$ are different. That is, calling by $\mu^t$ the law of $\pi_{Z(t)} M_\nu^t$, the following weak limits hold uniformly in $t$:  
$$ \lim_{x \rightarrow -\infty} \pi_x \mu^t \stackrel{\D}{=} \nu_{-}\,\, ,\,\,   \lim_{x \rightarrow +\infty} \pi_x \mu^t \stackrel{\D}{=} \nu_{+} \,,$$ 
 where $\nu_{-}$ has a density strictly smaller than the density of $\nu_{+}$. We will show that a second class particle initially located at the origin is a microscopic shock for both the two-class model $(\nu_\rho,\bar\nu_{\xi})$ and for the shock measure $\nu_{\rho,\lambda}$. The proof is based on a coupling between both regimes in such a way that the respective second-class particles stay close together. As a  consequence, we will also obtain the variance of the second-class particle in the two-class invariant regime. Let $(M^t_{\nu_\rho}, \bar{M}^t_{\xi'})_{t\geq 0} $ denote the two-class invariant process conditioned to have a second class particle at the origin, and let $Z_{\xi^{'}}(t)$ denote the position at time $t$ of the second-class particle which starts at the origin. 
 
\begin{theo}\label{theo:ZZ} 
There exists a coupling $(Z_{\nu'_{\rho,\lambda}} (t), Z_{\xi'} (t))$ such that 
\begin{equation}
 | Z_{\nu'_{\rho,\lambda}} (t) - Z_{\xi'} (t) | \leq J_t\mbox{ for all $t\geq 0$ }\,,
\end{equation}
where $(J_t)_{t\geq 0}$ is a family of identically distributed random variables with finite moments. Furthermore, $Z_{\xi'} (t)$ and $Z_{\nu'_{\rho,\lambda}}(t)$ are microscopic shocks.
\end{theo}

Together with Theorem \ref{Theo:Shock}, Theorem \ref{theo:ZZ} implies that:
\begin{coro}\label{coro:ZZ}
For all $u\in\RR$, 
$$\lim_{t\to\infty}\P\left(Z_{\xi'}(t)\leq \mu_2 t +(\sigma_2\sqrt{t}) u\right)=\Phi(u)\,.$$
For all $t \geq 0$,
\begin{equation}
\Var Z_{\xi'} (t) =  \sigma_2^2 t + O(t^{2/3})\,. \nonumber
\end{equation}

\end{coro}



\section{Proofs}

\subsection{First-class particles}
For $x<0$, we can decompose the flux $L_{\nu_\lambda}(x,t)$ as a difference  between two independent Poisson processes, by setting 
$$\Delta_\lambda(x,t):=L_{\nu_\lambda}(x,t) + \nu_{\lambda}((x,0])\,.$$
To see that $\Delta_\lambda(x,t)$ and $\nu_{\lambda}((x,0])$ are independent when $x<0$, note that $\Delta_\lambda(x,t)$ is a deterministic function of $\bP$ and $\nu_\lambda\mid_{(-\infty,x]}$ which are both independent of $\nu_{\lambda}((x,0])$. This becomes particularly clear when, as in \cite{CG1}, $(\Delta_\lambda(x,t)\,,\,t\geq 0)$ and $(-\nu_{\lambda}((x,0])\,,\,x\leq 0)$ are interpret as the west and south processes, respectively (associated to the space-time rectangle $(x,0]\times(0,t]$). Furthermore, the west process is a Poisson process of intensity $1/\lambda$ (for details, see \cite{CG1} and references there in). The same can be done to $L_{\nu'_\lambda}$ and this is the key to prove Theorem \ref{Theo:distX}.
   
\begin{demo}{\bf Proof of Theorem \ref{Theo:distX}.} 
Recall that $\nu'_\lambda:=\nu_\lambda+\delta_0$, and define for $t> 0$ and $x<0$
\begin{equation}\label{Eq:FluxAdd}
\Delta'_{\lambda}(x,t)=L_{\nu'_\lambda}(x,t) + \nu'_{\lambda}((x,0])\,.
\end{equation}
Since, for $x<0$, $\Delta'_\lambda(x,t)$ is a deterministic function of $\bP$ and $\nu'_\lambda\mid_{(-\infty,x]}$, which are both independent of $\nu'_{\lambda}((x,0])$, we obtain that $\Delta'_\lambda$ and $\nu'_\lambda$ are independent, with $\Delta'_\lambda$ being a Poisson processes of intensity $1/\lambda$ (see \cite{CG1}). By \eqref{Eq:FluxL},
\begin{equation}
\nonumber \left\{X_{\nu'_\lambda}(t) > x\right\} = \left\{L_{\nu'_\lambda} (x,t) < 0\right\} = \left\{\nu'_{\lambda}((x,0])> \Delta'_\lambda(x,t)\right\}\,,
\end{equation}
and hence, 
$$\P\left(X_{\nu'_\lambda}(t) > x\right)=\P\left(P_1(-x)+1>P_2(t)\right)\,,$$
where $P_1$ and $P_2$ are two independent Poisson processes with intensities $\lambda$ and $1/\lambda$, respectively, which completes the proof of \eqref{1dist}. 

\end{demo}


To compute the mean and the variance of $X_{\nu'_\lambda}(t)$ we use Lemma \ref{lemm:relation} below. An analogous result will be used to compute the mean and the variance of the second-class particle. Recall we labeled particles according to their initial positions, $X^z_{\nu_{\lambda}}(0)=z$ for $z\in\nu_{\lambda}$.   

\begin{lemm}\label{lemm:relation}
Define 
$$\chi(x,t):= \min \{z\in\nu\,:\, X_{\nu_\lambda}^z(t) > x\}\,.$$
Then for all $x\in\RR$ and $t\geq 0$ 
\begin{equation}\label{eq:psaidaX}
 \nu_\lambda(\chi(x,t)) =   L_{\nu_\lambda}(x,t) + 1
\end{equation}
Furthermore, for any measurable $f$ and $y< 0$ 
\begin{equation}\label{eq:econdX}
\E \left( f(\nu_\lambda(\chi(x,t))) \mid \chi(x,t)=y \right) = \E \left(f(\nu_\lambda(y))\right)\,,
\end{equation} 
and 
\begin{equation}\label{eq:relX}
\chi(x,t) \,{\buildrel \iD \over =}\, x - X_{\nu_\lambda}(t) \,.
\end{equation}
\end{lemm}

\begin{dem}
Notice that, for all $t\geq 0$ and $x\in\RR$ we have 
\begin{equation*}
 \nu_\lambda(\chi(x,t)) = \,\left\{\begin{array}{ll}  - M_{\nu_\lambda}^{t}( (x, X_{\nu_\lambda}(t)] ) + 1 & \mbox{if } x < X_{\nu_\lambda}(t)\,, \\
 1 & \mbox{if } x = X_{\nu_\lambda}(t) \,,\\
 M_{\nu_{\lambda}}^{t}( (X_{\nu_\lambda}(t), x] ) + 1 & \mbox{if } x > X_{\nu_\lambda}(t) \,.\end{array}\right.
\end{equation*}
By \eqref{Eq:FluxL}, we have \eqref{eq:psaidaX}.\\  

Now, if $\chi(x,t)=y$ and we fix $\bP$  and $\nu_\lambda\mid_{(-\infty,y]}$, any modification of $\nu_\lambda\mid_{(y,+\infty)}$, does not change the value of $\chi(x,t)$. To see this, let $z < y$ and observe that $X^z_{\nu_\lambda}(t)$ is independent of $\nu_\lambda\mid_{(y,+\infty)}$ since it is a deterministic function of $\bP$ and $\nu_\lambda\mid_{(-\infty,z]}$. Now, by definition 
$$\left\{\chi(x,t) = y \right\} = \left\{ \nu_\lambda(\{y\}) = 1\,,\,  X^y_{\nu_\lambda} (t) > x \right\}\cap \left\{X^z_{\nu_\lambda} (t) \leq x \,,\,  \forall z < y \text{ s.t. } \nu_\lambda(\{z\}) = 1 \right\}\,.$$ Therefore, $\left\{\chi (x,t) = y \right\}$ is independent of $\nu_\lambda\mid_{(y,+\infty)}$, and since $(y,0]\subseteq (y,+\infty)$ (recall that $y< 0$), it is also independent of $\nu_\lambda(y)$. Therefore, for a measurable $f$ and $y < 0$,  
\begin{equation}\nonumber
\E \left(f(\nu_\lambda(\chi(x,t))) \mid \chi(x,t) = y \right) = \E \left( f(\nu_\lambda(y))\right)\,,
\end{equation}
which gives \eqref{eq:econdX}. \\

Finally, from the definitions of $\chi_\nu$ and $X_{\nu}$ it follows that, 
\begin{equation}\label{eq:qui}
\left\{X_{\nu_\lambda}(t) \leq x\right\} = \left\{\chi(x,t) > 0\right\}\,
\end{equation}
for all $t\geq 0$ and $x\in\RR$. By translation invariance of both $\nu_\lambda$ and the two-dimensional Poisson process $\bP$, for any $h \in \RR$
\begin{eqnarray*}
\P\left(\chi(x,t) > h\right) &=& \P\left(\chi(x-h,t) > 0\right)\\ \nonumber
&=& \P\left(X_{\nu_\lambda}(t) \leq x - h \right) \\ \nonumber
&=& \P\left(x - X_{\nu_\lambda}(t) \geq h \right) \\ \nonumber
&=& \P\left(x - X_{\nu_\lambda}(t) > h \right)\,.
\end{eqnarray*}
This shows \eqref{eq:relX} and finishes the proof of the lemma.

\end{dem}


\begin{demo}{\bf Proof of Theorem \ref{Theo:cltX}.}
By \eqref{Eq:FluxL} and \eqref{eq:qui} we have $\left\{\chi(x,t) > 0\right\}\ \Rightarrow \left\{\L_{\nu_\lambda}(x,t)_{-} = 0\right\}$. Thus, we can focus on the case where $\chi(x,t)$ is non-positive. For $y< 0$
\begin{eqnarray*}
\E\left(L_{\nu_\lambda}(x,t)_{-} \mid \chi(x,t) = y \right) &=&  -\E\left(L_{\nu_\lambda}(x,t) \mid \chi(x,t) = y \right)
\end{eqnarray*}
and by \eqref{eq:psaidaX} and \eqref{eq:econdX} 
\begin{eqnarray*}
-\E\left(L_{\nu_\lambda}(x,t) \mid \chi(x,t) = y \right) &=&  -\E\left(\nu_\lambda(\chi(x,t)) - 1\mid \chi(x,t) = y\right) \\ \nonumber
&=& \lambda(-y) + 1
\end{eqnarray*}
which leads to
\begin{eqnarray*}
\E\left(L_{\nu_\lambda}(x,t)_{-}  \right) &=& \E\left(\E\left(L_{\nu_\lambda}(x,t)_{-} \mid \chi(x,t) \right)\right) \\ \nonumber
&=& \lambda\E\left(\chi(x,t)_{-} \right) + 1
\end{eqnarray*}

Now, \eqref{eq:relX} implies
\begin{equation}\label{eq:quineg}
\chi(x,t)_{-}  \,{\buildrel \iD \over =}\, \big(x - X_{\nu_\lambda}(t)\big)_{-}  \,
\end{equation}
then
\begin{eqnarray}\label{eq:meanrelation}
\E\left(L_{\nu_\lambda}(x,t)_{-}  \right) &=& \lambda\E\left([x- X_{\nu_{\lambda}}(t) ]_{-} \right) + 1
\end{eqnarray}

On the other hand, 
$$ 0\leq L_{\nu_\lambda}(x,t)_{{-}} + L_{\nu_\lambda}(x,t) = L_{\nu_\lambda}(x,t)_{+} \,,$$
and $L_{\nu_\lambda}(x,t)_{+}$ decreases to $0$ as $x\to-\infty$. By monotone convergence,
\begin{equation}
\lim_{x\to -\infty}\left\{\E\left(L_{\nu_\lambda}(x,t)_{-}\right) + \E\left(L_{\nu_\lambda}(x,t)\right) \right\} = 0 \nonumber\,.
\end{equation}
Analogously,
\begin{equation}
\lim_{x\to -\infty}\left\{\E\left([x- X_{\nu_{\lambda}}(t) ]_{-} \right) + \E\left( x- X_{\nu_{\lambda}}(t)  \right)\right\}=0 \nonumber\,.
\end{equation}
Together with \eqref{eq:meanrelation}, these three limits imply that
\begin{equation}
\lim_{x\to -\infty}\left\{\lambda\E\left(x- X_{\nu_{\lambda}}(t)  \right) - 1 - \E\left(L_{\nu_\lambda}(x,t)\right) \right\} = 0 \nonumber\,.
\end{equation} 
Now, by plugging 
$$\E\left(L_{\nu_\lambda}(x,t)\right)=\lambda x +\lambda^{-1}t\,,$$ 
in the above limit, one must get that 
$$\E X_{\nu_{\lambda}}(t)=-\lambda^{-2} t -1/\lambda\,.$$
By translation invariance, 
$$ X_{\nu'_\lambda} (t) \,{\buildrel \iD \over =}\,  X_{\nu_{\lambda}}(t) - X_{\nu_{\lambda}}(0)\,,$$
and since  $X_{\nu_{\lambda}}(0)$ is an exponential r.v. of rate $\lambda$, we obtain
$$\E X_{\nu'_\lambda} (t)= -\lambda^{-2} t\,.$$

For the second moment in \eqref{1mean&var}, we proceed as before and use Lemma \ref{lemm:relation} to deduce that, for $y \leq 0$ 
\begin{eqnarray*}
\E\left(L_{\nu_\lambda}(x,t)^{2}_{-} \mid \chi(x,t) = y \right) &=& \E\left((\nu_\lambda(\chi (x,t)) - 1)^{2} \mid \chi(x,t) = y \right) \\ \nonumber
&=& \lambda^2 y^2 + \lambda (-y) - 2\lambda y + 1
\end{eqnarray*}
Together with \eqref{eq:quineg} we obtain
\begin{eqnarray}\label{eq:varrelation}
\E\left(L_{\nu_\lambda}(x,t)^{2}_{-}\right) &=&  \E\left(\E\left(L_{\nu_\lambda}(x,t)^{2}_{-} \mid \chi(x,t) \right)\right)\\ \nonumber
&=& \lambda^2 \E\left(\chi(x,t)_{-}^2 \right) + 3\lambda \E\left(\chi(x,t)_{-} \right) + 1 \\ \nonumber
&=& \lambda^2 \E\left(([ x- X_{\nu_{\lambda}}(t) ]_{-})^2 \right) + 3\lambda \E\left([ x- X_{\nu_{\lambda}}(t) ]_{-} \right) + 1 \,.
\end{eqnarray}
Also
\begin{equation}
\lim_{x\to -\infty}\left\{\E\left(L_{\nu_\lambda}(x,t)^{2}_{-}\right) - \E\left(L_{\nu_\lambda}(x,t)^2\right) \right\} = 0 \nonumber\,,
\end{equation}
and 
\begin{equation}
\lim_{x\to -\infty}\left\{\E\left(([ x- X_{\nu_{\lambda}}(t) ]_{-})^2 \right)- \E\left( (x- X_{\nu_{\lambda}}(t))^2  \right)\right\}  = 0 \nonumber\,.
\end{equation}
Together with \eqref{eq:varrelation}, these yield to
\begin{equation}
\lim_{x\to -\infty}\left\{\lambda^2\E\left( (x- X_{\nu_{\lambda}}(t))^2  \right)  - 3\lambda \E\left( x- X_{\nu_{\lambda}}(t)  \right) + 1 - \E\left(L_{\nu_\lambda}(x,t)^2\right)  \right\}= 0 \nonumber
\end{equation}
Since for $x \leq 0$ and $t \geq 0$, $L_{\nu_{\lambda}}(x,t)$ is the difference between  two independent Poisson processes of rates $\frac{t}{\lambda}$ and $-\lambda x$, respectively, we obtain
$$\E\left(L_{\nu_\lambda}(x,t)^2\right) = \lambda^2 x^2 -\lambda x  + 2xt  + \frac{t^2}{\lambda^2}  + \frac{t}{\lambda}\,,$$
and one must have that 
\begin{equation}
\E\left( X_{\nu_{\lambda}}(t)^2 \right)  = \frac{t^2}{\lambda^4} + \frac{4t}{\lambda^3} + \frac{2}{\lambda^2}  \nonumber\,,
\end{equation}
and hence,
\begin{equation}
\Var X_{\nu_{\lambda}}(t) = \frac{2t}{\lambda^3} + \frac{1}{\lambda^2}   \nonumber\,.
\end{equation}
Since $X_{\nu_{\lambda}}(t)- X_{\nu_{\lambda}}(0)$ is a deterministic function of $\bP$ and $\nu_{\lambda}\mid_{(-\infty, X_{\nu_{\lambda}}(0)]}$, it is independent of $\nu_{\lambda}\mid_{(X_{\nu_{\lambda}}(0), +\infty)}$. In particular, $X_{\nu_{\lambda}}(t) - X_{\nu_{\lambda}}(0)$ and $X_{\nu_{\lambda}}(0)$ are independent, which leads to
\begin{equation}
\Var X_{\nu'_\lambda}(t)=\Var\left(X_{\nu_{\lambda}}(t) - X_{\nu_{\lambda}}(0)\right) = \frac{2t}{\lambda^3}  \nonumber\,.
\end{equation}

\end{demo}



\subsection{Asymptotics for the variance of flux of second-class particles.} In this section we prove Theorem \ref{Theo:2fluxVar}. The key ingredient is a more general formulation of the exit point formula, previously proved in \cite{CG2}, that relates the covariance 
between the initial profile and the flux (of first class particles) with the expected position of the exit point $Y_{\nu_\lambda}$:
\begin{equation*}
\Cov\left(L_{\nu_\lambda}(x,t),\nu_\lambda(x)\right)=\E Y^s_{\nu_\lambda}(x,t)_+\,.
\end{equation*}
   
\begin{lemm}\label{exitpoint}
Fix $\rho\leq \lambda$ and define $\nu_\rho$ by selecting points of $\nu_\lambda$ independently with probability $\rho\lambda^{-1}$, so that $\nu_\lambda$ can be seen as the sum of two independent Poisson process $\nu_\rho$ and $\nu_{\lambda-\rho}$. If $y \geq 0$ then 
$$\Cov\left(L_{\nu_\lambda}(x,t),\nu_\rho(y)\right)=\rho\E\left(Y^s_{\nu_\lambda}(x,t)_+\wedge y\right)\,,$$
while if $y<0$ then
$$\Cov\left(L_{\nu_\lambda}(x,t),\nu_\rho(y)\right)=\rho\E\left(Y^i_{\nu_\lambda}(x,t)_-\wedge (-y)\right)\,.$$
\end{lemm}

\begin{dem}
To prove the first formula we add a Poisson process on the interval $(0,y]$ of intensity $\epsilon$ to $\nu_\rho$, and denote the flux obtained from this new initial profile by $L^\epsilon_\lambda$. We are seeing the flux as a function of $\nu_{\rho+\epsilon}$ and $\nu_{\lambda-\rho}$. Then we compute the derivative $f_+'(\epsilon) = \partial_\epsilon\E L_{\nu_\lambda}^{+,\epsilon}(x,t)$ in two different ways. First we condition on the value of $\nu_{\rho+\epsilon}(y)$ and denote $a_{n}:=\E\big(L_{\nu_\lambda}^{+,\epsilon}|\nu_{\rho+\epsilon}(y) = n \big)$. Since   $a_{n}$ does not depend on $\epsilon$ anymore, we have that
\begin{align*}
f_+'(0) &= \partial_\epsilon|_{\epsilon = 0}\sum_{n=0}^{\infty}\frac{((\rho +\epsilon)t)^n}{n!}e^{-t(\rho +\epsilon)}a_{n} \\
 &= \frac{1}{\rho}\big(\sum_{n=0}^{\infty}\frac{(\rho t)^n}{n!}e^{-t\rho} a_{n} n\big)  - t\sum_{n=0}^{\infty}\frac{(\rho t)^n}{n!}e^{-t\rho} a_{n}\\
 &= \frac{1}{\rho}\Cov\left(L_{\nu_\lambda}(x,t),\nu_\rho(y)\right)
\end{align*}

The variational representation \eqref{last-evo} of the flux tells us that, if one adds only one point to $(0,y]$, which is the order $\epsilon$ term in the Taylor expansion of $f(\epsilon)$, then the flux will increase by one if and only if this point belongs to $(0,Y^s_{\nu_\lambda}(x,t)_+]$. Thus, 
$$f_+(\epsilon)=f_+(0)+\epsilon\E\left(Y^s_{\nu_\lambda}(x,t)_+\wedge y\right) + O(\epsilon^2)\,,$$
which yields to 
$$\Cov\left(L_{\nu_\lambda}(x,t),\nu_\rho(y)\right)=\rho\E\left(Y^s_{\nu_\lambda}(x,t)_+\wedge y\right)\,.$$

For $y<0$, we add to $\nu_\rho$ a Poisson process on the interval $(y,0]$ of intensity $\epsilon$ to get a new profile. For this new $f_-$, by conditioning on the values of $-\nu_{\rho+\epsilon}(y)$ (the number of points in $(y,0]$), we get   
$$f'_-(0)=\frac{1}{\rho}\Cov\left(L_{\nu_\lambda}(x,t),-\nu_\rho(y)\right)=-\frac{1}{\rho}\Cov\left(L_{\nu_\lambda}(x,t),\nu_\rho(y)\right)\,.$$
On the other hand, if one add only one point in $(y,0]$, then the flux will decrease by one if and only if this point belongs to $(-Y^i_{\nu_\lambda}(x,t)^-,0]$ (again using \eqref{last-evo}) and so 
$$f_-(\epsilon)=f_-(0)-\epsilon\E\left(Y^i_{\nu_\lambda}(x,t)_-\wedge (-y)\right) + O(\epsilon^2)\,,$$
which proves the second covariance formula.
   
\end{dem}


\begin{demo}{\bf Proof of Theorem \ref{Theo:2fluxVar}.}
Recall that 
$$\Gamma_{\rho,\lambda}(x,t):=\nu_\lambda(x-\lambda^{-2}t)-\nu_\rho(x-\rho^{-2}t)+2(\lambda^{-1}-\rho^{-1})t\,.$$
Then (we are omitting the dependence on $(x,t)$), 
\begin{equation*}
\nonumber \Var \bar L_{\rho,\lambda}=\Var\Gamma_{\rho,\lambda}+\Var\left(\bar L_{\rho,\lambda}-\Gamma_{\rho,\lambda}\right) +2\Cov\left(\bar L_{\rho,\lambda}-\Gamma_{\rho,\lambda}\,,\,\Gamma_{\rho,\lambda}\right)\,.
\end{equation*}
By translation invariance (see \eqref{Theo:Flux2}),
$$\Var\left(\bar L_{\rho,\lambda}(x,t)-\Gamma_{\rho,\lambda}(x,t)\right)=\E\left(\left\{\bar L_{\rho,\lambda}(x,t)-\Gamma_{\rho,\lambda}(x,t)\right\}^2\right) =R_1(t)\,.$$
On the other hand, 
\begin{eqnarray*}
\Var\Gamma_{\rho,\lambda}(x,t)&=&\Var\left(\nu_\lambda(x-\lambda^{-2}t)-\nu_\rho(x-\rho^{-2}t)\right)\\
&=&\Var\bar \nu_{\rho,\lambda}(x-\rho^{-2}t)+\Var\nu_\lambda(\rho^{-2}t-\lambda^{-2}t)\\
&=&\left(\lambda-\rho\right)\left(x-\rho^{-2}t\right)+\lambda\left(\rho^{-2}t-\lambda^{-2}t\right)\\
&=&\left(\lambda-\rho\right)\left(x+\frac{t}{\lambda\rho}\right)\,.
\end{eqnarray*}
Hence, we need to show that 
$$\Cov\left(\bar L_{\rho,\lambda}-\Gamma_{\rho,\lambda}\,,\,\Gamma_{\rho,\lambda}\right)=R_2+R_3\,.$$
To do so, let
$$\tilde L_{\nu_\lambda}(x,t):=L_{\nu_\lambda}(x,t)-\nu_\lambda(x-\lambda^{-2}t)\,,$$
and 
$$\tilde L_{\nu_\rho}(x,t):=L_{\nu_\rho}(x,t)-\nu_\rho(x-\rho^{-2}t)\,.$$
Then
$$\Cov\left(\bar L_{\rho,\lambda}-\Gamma_{\rho,\lambda}\,,\,\Gamma_{\rho,\lambda}\right)=\Cov\left(\tilde L_{\nu_\rho},\nu_\rho\right)-\Cov\left(\tilde L_{\nu_\rho},\nu_\lambda\right)+\Cov\left(\tilde L_{\nu_\lambda},\nu_\lambda\right)-\Cov\left(\tilde L_{\nu_\lambda},\nu_\rho\right)\,.$$

Now, by using translation invariance together with the second equation in the statement of Lemma \ref{exitpoint} (since $\rho^{-2}t-x\leq 0$), 
\begin{eqnarray*}
\Cov\left(\tilde L_{\nu_\rho}(x,t),\nu_\rho(x-\rho^{-2}t)\right)&=&-\Cov\left(L_{\nu_\rho}(\rho^{-2}t,t),\nu_\rho(\rho^{-2}t-x)\right)\\
&=&-\rho\E\left(Y^i_{\nu_\rho}(\rho^{-2}t,t)_-\wedge(x-
\rho^{-2}t)\right)\,.
\end{eqnarray*}
Analogously (recall that $\nu_\lambda=\nu_\rho+\bar\nu_{\lambda-\rho}$),
\begin{eqnarray*}
\Cov\left(\tilde L_{\nu_\rho}(x,t),\nu_\lambda(x-\lambda^{-2}t)\right)&=&\Cov\left(L_{\nu_\rho}(\rho^{-2}t,t),\nu_\rho(\rho^{-2}t-\lambda^{-2}t)\right)\\
&-&\Cov\left(L_{\nu_\rho}(\rho^{-2}t,t),\nu_\rho(\rho^{-2}t-x)\right)\\
&=&\rho\E\left(Y^s_{\nu_\rho}(\rho^{-2}t,t)_+\wedge(\rho^{-2}t-
\lambda^{-2}t)\right)\\
&-&\rho\E\left(Y^i_{\nu_\rho}(\rho^{-2}t,t)_-\wedge(x-
\rho^{-2}t)\right)\,,
\end{eqnarray*}
which shows that
$$\Cov\left(\tilde L_{\nu_\rho},\nu_\rho\right)-\Cov\left(\tilde L_{\nu_\rho},\nu_\lambda\right)=R_2\,.$$
By repeating the same reasoning as before,
\begin{equation*}
\Cov\left(\tilde L_{\nu_\lambda}(x,t),\nu_\lambda(x-\lambda^{-2}t)\right)=-\lambda\E\left(Y^i_{\nu_\lambda}(\lambda^{-2}t,t)_-\wedge(x-
\lambda^{-2}t)\right)\,,
\end{equation*}
and
\begin{eqnarray*}
\Cov\left(\tilde L_{\nu_\lambda}(x,t),\nu_\rho(x-\rho^{-2}t)\right)&=&\Cov\left(L_{\nu_\lambda}(\lambda^{-2}t,t),\nu_\rho(\lambda^{-2}t-\rho^{-2}t)\right)\\
&-&\Cov\left(L_{\nu_\lambda}(\lambda^{-2}t,t),\nu_\rho(\lambda^{-2}t-x)\right)\\
&=&\rho\E\left(Y^i_{\nu_\lambda}(\lambda^{-2}t,t)_-\wedge(\rho^{-2}t-\lambda^{-2}t)\right)\\
&-&\rho\E\left(Y^i_{\nu_\lambda}(\lambda^{-2}t,t)_-\wedge(x-\lambda^{-2}t)\right)\,,
\end{eqnarray*}
and thus
$$\Cov\left(\tilde L_{\nu_\lambda},\nu_\lambda\right)-\Cov\left(\tilde L_{\nu_\lambda},\nu_\rho\right)=R_3\,.$$

\end{demo}

\subsection{Second-class particles}
The proof of Theorem \ref{Theo:Shock} is similar to what we have done so far. Denote by $Z_{\bar\nu_{\lambda-\rho}}(t)$ the position at time $t$ of the second class particle that was at the first particle $a\in\bar\nu_{\lambda-\rho}:=\nu_\lambda-\nu_\rho$ (a Poisson process of rate $\lambda - \rho$ independent of $\nu_\rho$) such that $a\leq 0$. The key is again to relate the position of the tagged particle with the flux of particles: for all $t\geq 0$ and $x\in\RR$: 
\begin{equation}\label{Eq:2FluxL}
 \bar L_{\rho,\lambda}(x,t) = \,\left\{\begin{array}{ll}  - \bar M_{\bar\nu_{\lambda-\rho}}^{t}( (x, Z_{\bar\nu_{\lambda-\rho}}(t)] ) & \mbox{if } x < Z_{\bar\nu_{\lambda-\rho}}(t)\,, \\
 0 & \mbox{if } x = Z_{\bar\nu_{\lambda-\rho}}(t)\,, \\
 \bar M_{\bar\nu_{\lambda-\rho}}^{t}( (Z_{\bar\nu_{\lambda-\rho}}(t), x] ) & \mbox{if } x > Z_{\bar\nu_{\lambda-\rho}}(t)\,. \end{array}\right.
\end{equation}
Thus, $Z_{\bar\nu'_{\lambda-\rho}(t)}=Z_{\nu'_{\rho,\lambda}}(t)$ denotes the position at time $t$ of the tagged second class particle whose initial position is $0\in\bar\nu'_{\lambda-\rho}:=\bar\nu_{\lambda-\rho}+\delta_0$.
 
\begin{demo}{\bf Proof of Theorem \ref{Theo:CLTShock}.}
By \eqref{Eq:2FluxL},
$$\P\left(Z_{\nu'_{\rho,\lambda}}(t)\leq x\right)=\P\left(L_{\nu_\rho}(x,t) \leq L_{\nu'_\lambda}(x,t)\right)\,.$$
Now, for $t>0$, define 
$$w_u(t)=\mu_2 t + (\sigma_2 \sqrt{t})u>0\,.$$
Then
$$\P\left(Z_{\nu'_{\rho,\lambda}}(t)\leq \mu_2 t +(\sigma_2\sqrt{t}) u\right)=\P\left(L_{\nu_\rho}(w_u(t),t) \leq L_{\nu'_\lambda}(w_u(t),t)\right)\,.$$
Therefore, we only need to prove that 
\begin{equation}\label{eq:2clt}
\lim_{t\to\infty}\P\left(L_{\nu_\rho}(w_u(t),t) \leq L_{\nu'_\lambda}(w_u(t),t)\right) =\Phi(u)\,.
\end{equation}

To prove \eqref{eq:2clt}, let $\eta:=\sqrt{\frac{\lambda-\rho}{\lambda\rho}}$, and define
$$\Delta_{u,\lambda}(t):=\frac{L_{\nu'_\lambda}(w_u(t),t)-\left(\lambda w_u(t)+\frac{t}{\lambda}\right)}{\eta\sqrt{t}}\,\,\mbox{ and }\,\,\Delta_{u,\rho}(t):=\frac{L_{\nu_\rho}(w_u(t),t)-\left(\rho w_u(t)+\frac{t}{\rho}\right)}{\eta\sqrt{t}}\,.$$
Since 
$$(\lambda - \rho)w_u(t) + \left(\lambda^{-1} - \rho^{-1}\right)t=u \eta \sqrt{2}\sqrt{t}\,,$$
then
\begin{equation}\label{lem:comparation}
\left\{L_{\nu_\rho}(w_u(t),t)\leq L_{\nu'_\lambda}(w_u(t),t)\right\}=\left\{\frac{\Delta_{u,\rho}(t)-\Delta_{u,\lambda}(t)}{\sqrt{2}}\leq u\right\}\,.
\end{equation}
Define 
\begin{equation*}
A_{u,\lambda} (t):= \frac{L_{\nu'_\lambda}(w_u(t),t)-\left[\nu_\lambda(w_u(t) - t/\lambda^2)+\frac{2}{\lambda}t\right]}{\eta\sqrt{2}\sqrt{t}}\,,
\end{equation*} 
\begin{equation*}
A_{u,\rho} (t):= \frac{L_{\nu_\rho}(w_u(t),t)-\left[\nu_\rho(w_u(t) - t/\rho^2)+\frac{2}{\rho}t\right]}{\eta\sqrt{2}\sqrt{t}}\,, 
\end{equation*}
and, for $\psi=\lambda,\rho$,
\begin{equation}\nonumber
B_{u,\psi} := \frac{\left[\nu_\psi(w_u(t) - t/\psi^2)+\frac{2}{\psi}t\right] - \E\left(\nu_\psi(w_u(t) - t/\psi^2)+\frac{2}{\psi}t\right)}{\eta\sqrt{2}\sqrt{t}} \,.
\end{equation}
Then
$$\frac{\Delta_{u,\rho}(t)-\Delta_{u,\lambda}(t)}{\sqrt{2}}  = \left[A_{u,\rho}(t) - A_{u,\lambda}(t)\right] - \left[B_{u,\lambda}(t) - B_{u,\rho}(t)\right]\,.$$
By \eqref{1classInf} for $\psi=\lambda,\rho$ (and fixed $u\in\RR$)
$$\lim_{t\to\infty}\E(A_{u,\psi}^2) = 0\,,$$
and therefore, 
\begin{equation}\label{eq:influence}
\lim_{t\to\infty}\E\left(\Big\{A_{u,\rho}(t) - A_{u,\lambda}(t)\Big\}^2\right)=0\,.
\end{equation}
(Since $\nu'_\lambda=\nu_\lambda+\delta_0$, \eqref{1classInf} also holds for $L_{\nu'_\lambda}$.) By the central limit theorem for Poisson process, for $\psi=\lambda,\rho$ (and fixed $u\in\RR$), $\sqrt{2}B_{u,\psi}(t)$ converges, in distribution, to a standard Gaussian random variable. Notice also that, for sufficiently large $t>0$, 
$$w_u(t) - \frac{t}{\lambda^2} > 0\,\,\mbox{ and }\,\,w_u(t) - \frac{t}{\rho^2} < 0\,.$$
This implies that $B_{u,\rho}(t)$ and $B_{u,\lambda}(t)$ are independent, for sufficiently large $t>0$. Hence, $B_{u,\lambda}(t) - B_{u,\rho}(t)$ converges, in distribution, to a standard Gaussian random variable. Together with \eqref{lem:comparation} and \eqref{eq:influence}, this finishes the proof of \eqref{eq:2clt}.

\end{demo}

As before, to compute first and second moments of the tagged second class particle we need an extra lemma. Recall we labeled particles according to their initial positions, $Z^y_{\bar\nu_{\lambda-\rho}}(0)=y$ for $y\in\bar\nu_{\lambda-\rho}$. We denote by $\bar Z_{\bar\nu_{\lambda-\rho}} (t)$ the position at time $t$ of the first second-class particle strictly to the right of the origin.   
\begin{lemm}\label{lemm:relation2}
Define  
$$\bar\chi(x,t) := \max \{y\in \bar\nu_{\lambda-\rho} \,:\, Z^{y}_{\bar\nu_{\lambda-\rho}}(t) \leq x\}\,.$$ 
Then
\begin{equation}\label{eq:psaidaZ}
\bar L_{\rho,\lambda}(x,t) = \bar\nu_{\lambda-\rho}\left(\bar\chi(x,t)\right)\,.
\end{equation}
Furthermore, for any measurable $f$ and $y > 0$, 
\begin{equation}\label{eq:econdZ}
\E \left( f(\bar\nu_{\lambda-\rho}(\bar\chi(x,t) )) \mid\bar\chi(x,t) = y \right)=\E \left(f(\bar\nu_{\lambda-\rho}(y)+1)\right)\,,
\end{equation} 
and
\begin{equation}\label{eq:relZ}
\bar\chi(x,t) \,\,{\buildrel \iD \over =}\,\, x-\bar Z_{\bar\nu_{\lambda-\rho}} (t)\,.
\end{equation}
 \end{lemm}

\begin{dem}
The proof is basically the same of Lemma \ref{lemm:relation}. For all $t\geq 0$ and $x\in\RR$ we have 
\begin{equation*}
 \nu(\bar\chi(x,t)) = \,\left\{\begin{array}{ll}  - \bar M_{\bar\nu_{\lambda-\rho}}^{t}( (x, Z_{\bar\nu_{\lambda-\rho}}(t)] ) & \mbox{if } x < Z_{\bar\nu_{\lambda-\rho}}(t) \\
 0 & \mbox{if } x = Z_{\bar\nu_{\lambda-\rho}}(t) \\
 \bar M_{\bar\nu_{\lambda-\rho}}^{t}( (Z_{\bar\nu_{\lambda-\rho}}(t), x] ) & \mbox{if } x > Z_{\bar\nu_{\lambda-\rho}}(t) \end{array}\right.
\end{equation*}
By \eqref{Eq:2FluxL}, this implies \eqref{eq:psaidaZ}.\\ 

If $\bar\chi(x,t)=y$ and we fix $\bP$, $\nu_\rho$ and $\bar\nu_{\lambda-\rho}\mid_{[y,+\infty)}$,  any modification of $\bar\nu_{\lambda-\rho}\mid_{(-\infty, y)} $ does not change the value of $\bar\chi(x,t)$. Indeed, the trajectory of a tagged second-class particle with initial position $z > y$ is by definition a deterministic function of $\bP$, $\nu_\rho$  and $\bar\nu_{\lambda-\rho}\mid_{[z,+\infty)}$. Since $\bar\nu_{\lambda-\rho}$ is a Poisson process, $\bar\nu_{\lambda-\rho}\mid_{[z,+\infty)}$ is independent of $\bar\nu_{\lambda-\rho}\mid_{(-\infty, y)} $ (and independent of $\bP$ and $\nu_\rho$). Now, by definition
$$\left\{\bar\chi(x,t) = y \right\} = \left\{ \bar\nu_{\lambda-\rho}(\{y\}) = 1,  Z^{y}_{\bar\nu_{\lambda-\rho}} (t) \leq x \right\}\cap \left\{Z^{z}_{\bar\nu_{\lambda-\rho}} (t) > x ,  \forall z > y \text{ s.t. } \bar\nu_{\lambda-\rho}(\{z\}) = 1 \right\} \,.$$
Therefore, $\left\{\bar\chi (x,t) = y \right\}$ is independent of $\bar\nu_{\lambda-\rho}\mid_{(-\infty, y)} $, and since $(0,y)\subseteq (-\infty,y)$ (recall that $y>0$), it is also independent of $\nu_{\lambda}((0,y))$. Thus, for any measurable $f$ and $y > 0$,  
\begin{equation}\nonumber
\E \left( f(\bar\nu_{\lambda-\rho}(\bar\chi (x,t) )) \mid\bar\chi (x,t) = y \right) = \E \left( f(\bar\nu_{\lambda-\rho}(y)+1)\right)\,.
\end{equation}
(Notice that $\bar\chi (x,t) = y$ means that we have a particle at $y$.)
  
Finally, from the definitions of $\chi_{\rho,\lambda}$ and $\bar Z_{\bar\nu_{\lambda-\rho}}$, it follows that 
\begin{equation}\label{eq:quiZ} 
\left\{\bar Z_{\bar\nu_{\lambda-\rho}} (t) \leq x\right\} = \left\{\bar\chi (x,t) > 0\right\}\,
\end{equation}
for all $t\geq 0$ and $x\in\RR$. By translation invariance of $\nu_\rho$, $\bar\nu_{\lambda - \rho}$ and of the two-dimensional Poisson process $\bP$, for any constant $h \in \RR$
\begin{eqnarray*}
\P\left(\bar\chi(x,t) > h\right) &=& \P\left(\bar\chi (x-h,t) > 0\right)\\ \nonumber
&=& \P\left(\bar Z_{\bar\nu_{\lambda-\rho}} (t) \leq x - h \right) \\ \nonumber
&=& \P\left(x - \bar Z_{\bar\nu_{\lambda-\rho}} (t) \geq h \right) \\ \nonumber
&=& \P\left(x - \bar Z_{\bar\nu_{\lambda-\rho}} (t) > h \right)\,.
\end{eqnarray*}
This shows \eqref{eq:relZ} and finishes the proof of the lemma.
  
\end{dem}

\begin{demo}{\bf Proof of Theorem \ref{Theo:Shock}.} By \eqref{eq:quiZ} and the definition of $\bar L_{\rho,\lambda}(x,t)$ we have that $\left\{\bar\chi(x,t) < 0\right\}\ \Rightarrow \left\{\bar L_{\rho,\lambda}(x,t)_{+} = 0\right\}$. Thus, we can focus on the case where $\bar\chi(x,t)$ is non-negative. Fix $y> 0$, by \eqref{eq:psaidaZ} and \eqref{eq:econdZ} 
\begin{eqnarray*}
\E\left(\bar L_{\rho,\lambda}(x,t)_{+} \mid \bar\chi(x,t) = y \right) &=&  \E\left(\bar L_{\rho,\lambda}(x,t) \mid \bar\chi(x,t) = y \right) \\ \nonumber
&=&  \E\left(\bar\nu_{\lambda-\rho}(\bar\chi(x,t)) \mid \bar\chi(x,t) = y\right) \\ \nonumber
&=& (\lambda - \rho)y + 1
\end{eqnarray*}
which leads to
\begin{eqnarray*}
\E\left(\bar L_{\rho,\lambda}(x,t)_{+}  \right) &=& \E\left(\E\left(\bar L_{\rho,\lambda}(x,t)_{+} \mid \bar\chi(x,t) \right)\right) \\ \nonumber
&=& (\lambda - \rho)\E\left(\bar\chi(x,t)_{+} \right) + 1
\end{eqnarray*}

Now, \eqref{eq:relZ} implies
\begin{equation}\label{eq:quiZneg}
\bar\chi(x,t)_{+}  \,{\buildrel \iD \over =}\, \big(x - \bar Z_{\bar\nu_{\lambda-\rho}} (t)\big)_{+}  \,,
\end{equation}
and then
\begin{eqnarray}\nonumber
\E\left(\bar L_{\rho,\lambda}(x,t)_{+} \right) &=& (\lambda - \rho)\E\left([x- \bar Z_{\bar\nu_{\lambda-\rho}} (t) ]_{+} \right) + 1\,.
\end{eqnarray}

By monotone convergence,
\begin{equation}\nonumber
\lim_{x\to \infty}\left\{\E\left(\bar L_{\rho,\lambda}(x,t)_{+}\right) - \E\left(\bar L_{\rho,\lambda}(x,t)\right)  \right\}= 0 \,,
\end{equation}
and
\begin{equation}
\lim_{x\to \infty}\left\{\E\left([ x- \bar Z_{\bar\nu_{\lambda-\rho}} (t)]_{+} \right) - \E\left(x- \bar Z_{\bar\nu_{\lambda-\rho}} (t)  \right)  \right\}= 0\,, \nonumber
\end{equation}
which imply 
\begin{equation}
\lim_{x\to \infty}\left\{(\lambda -\rho)\E\left( x- \bar Z_{\bar\nu_{\lambda-\rho}} (t)  \right) + 1 - \E\left(\bar L_{\rho,\lambda}(x,t)\right) \right\} = 0\,. \nonumber
\end{equation}
Since
$$\E\left(\bar L_{\rho,\lambda}(x,t)\right) = \E\left(L_{\nu_\lambda}(x,t)\right) - \E\left(L_{\nu_\rho}(x,t)\right)= \left(\lambda - \rho\right)\left(x - \frac{t}{\lambda\rho} \right)\,,$$
then
\begin{equation}
\lim_{x\to \infty}\left\{(\lambda -\rho)\E\left(x -  \bar Z_{\bar\nu_{\lambda-\rho}} (t) \right) + 1 - (\lambda -\rho)(x - \frac{t}{\lambda\rho})\right\}  = 0 \,,\nonumber
\end{equation}
and thus
\begin{equation*}
\E\left( \bar Z_{\bar\nu_{\lambda-\rho}} (t)\right)  = \frac{t}{\lambda\rho} + \frac{1}{\lambda-\rho} \,. \nonumber
\end{equation*}
By translation invariance,
 $$Z_{\nu'_{\rho,\lambda}}(t)=Z_{\bar\nu'_{\lambda-\rho}}(t)= \,\,{\buildrel \iD \over =}\,\, \bar Z_{\bar\nu_{\lambda-\rho}} (t) - \bar Z_{\bar\nu_{\lambda-\rho}} (0)\,,$$
and since $\bar Z_{\bar\nu_{\lambda-\rho}} (0)$ is an exponential r.v. of rate $\lambda - \rho$, we obtain
\begin{equation}
\E\left( Z_{\nu'_{\rho,\lambda}}(t)\right)  = \frac{t}{\lambda\rho}  \,. \nonumber
\end{equation}\\


For the second moment we proceed as before and use Lemma \ref{lemm:relation2}: for $y > 0$, 
\begin{eqnarray*}
\E\left(\bar L_{\rho,\lambda}(x,t)^{2}_{+} \mid \bar\chi(x,t) = y \right) &=& \E\left(\bar\nu_{\lambda-\rho}(\bar\chi (x,t))^{2} \mid \bar\chi(x,t) = y \right) \\ \nonumber
&=& (\lambda -\rho)^2 y^2 + (\lambda -\rho)y + 2(\lambda -\rho)y + 1\,,
\end{eqnarray*}
and hence, 
\begin{eqnarray*}
\E\left(\bar L_{\rho,\lambda}(x,t)^{2}_{+}\right) &=& (\lambda -\rho)^2 \E\left(\bar\chi(x,t)_{+}^2 \right) + 3(\lambda -\rho)\E\left(\bar\chi(x,t)_{+} \right) + 1 \\ \nonumber
&=& (\lambda -\rho)^2 \E\left(([ x- \bar Z_{\bar\nu_{\lambda-\rho}} (t) ]_{+})^2 \right) + 3(\lambda -\rho)\E\left([ x- \bar Z_{\bar\nu_{\lambda-\rho}} (t) ]_{+} \right) + 1 \,.
\end{eqnarray*}

By monotone convergence, 
\begin{equation}
\lim_{x\to \infty}\left\{\E\left(\bar L_{\rho,\lambda}(x,t)^{2}_{+}\right) - \E\left(\bar L_{\rho,\lambda}(x,t)^2\right) \right\} = 0\,, \nonumber
\end{equation}
and
\begin{equation}
\lim_{x\to \infty}\left\{\E\left(([ x- \bar Z_{\bar\nu_{\lambda-\rho}} (t) ]_{+})^2 \right) - \E\left((x- \bar Z_{\bar\nu_{\lambda-\rho}} (t))^2 \right)\right\} = 0 \,,\nonumber
\end{equation}
which implies that
\begin{equation}
\lim_{x\to \infty}\left\{(\lambda - \rho)^2\E\left((x- \bar Z_{\bar\nu_{\lambda-\rho}} (t))^2  \right)  + 3(\lambda - \rho)\E\left( x- \bar Z_{\bar\nu_{\lambda-\rho}} (t)  \right) + 1  - \E\left(\bar L_{\rho,\lambda}(x,t)^2\right)\right\}  = 0\,. \nonumber
\end{equation}\\
By using that $\E\bar Z_{\bar\nu_{\lambda-\rho}} (t)=t(\lambda\rho)^{-1} + (\lambda-\rho)^{-1}$,
$$(\lambda - \rho)^2\E\left((x- \bar Z_{\bar\nu_{\lambda-\rho}} (t))^2  \right)= (\lambda - \rho)^2 \left(\Var \bar Z_{\bar\nu_{\lambda-\rho}} (t) + \left(x - \frac{t}{\lambda\rho} -\frac{1}{\lambda-\rho}\right)^2 \right)\,,$$
$$\E\left(\bar L_{\rho,\lambda}(x,t)^2\right)=  \Var\bar L_{\rho,\lambda}(x,t) + (\lambda-\rho)^2\left(x - \frac{t}{\lambda\rho}\right)^2\,,$$
(together with the previous limit) one finally gets 
\begin{equation*}
(\lambda - \rho)^2\Var \bar Z_{\bar\nu_{\lambda-\rho}} (t) = \lim_{x\to \infty}  \left\{\Var \bar L_{\rho,\lambda}(x,t) - (\lambda - \rho) \left(x - \frac{t}{\lambda\rho} \right) + 1 \right\}\,.   
\end{equation*}
By \eqref{2fluxAsymp}, we have 
\begin{eqnarray}
(\lambda - \rho)^2\Var \bar Z_{\bar\nu_{\lambda-\rho}} (t)  &=& \lim_{x\to \infty}\left\{  \Var\bar L_{\rho,\lambda}(x,t) - (\lambda - \rho) \left(x - \frac{t}{\lambda\rho} \right) + 1 \right\} \nonumber \\
&=& 2\big( \lambda - \rho\big)\frac{t}{\lambda \rho} + 1 +\lim_{x\to \infty} \left\{ \Var\bar L_{\rho,\lambda}(x,t)-\left( \lambda - \rho\right)\left(x + \frac{t}{\lambda \rho} \right)\right\}\nonumber \\
&=& 2\big( \lambda - \rho\big)\frac{t}{\lambda \rho} + 1 + O(t^{2/3})\,.\nonumber
\end{eqnarray}
Therefore
\begin{equation*}
\Var \bar Z_{\bar\nu_{\lambda-\rho}} (t) = \frac{2t}{(\lambda - \rho)\lambda \rho} + \frac{1}{(\lambda - \rho)^2} + O(t^{2/3}) \,.   
\end{equation*}
Now,
$$ Z_{\nu'_{\rho,\lambda}}(t) \,{\buildrel \iD \over =}\,  \bar Z_{\bar\nu_{\lambda-\rho}} (t) - \bar Z_{\bar\nu_{\lambda-\rho}} (0)\,.$$
Since $\bar Z_{\bar\nu_{\lambda-\rho}} (t) - \bar Z_{\bar\nu_{\lambda-\rho}} (0)$ is a deterministic function of $\bP$, $\nu_\rho$ and $\bar\nu_{\lambda-\rho}\mid_{[\bar Z_{\bar\nu_{\lambda-\rho}} (0), +\infty)}$, it is independent of $\bar\nu_{\lambda-\rho}\mid_{(-\infty, \bar Z_{\bar\nu_{\lambda-\rho}} (0))}$. In particular, $\bar Z_{\bar\nu_{\lambda-\rho}} (t) - \bar Z_{\bar\nu_{\lambda-\rho}} (0)$ and $\bar Z_{\bar\nu_{\lambda-\rho}} (0)$ are independent, which leads to
\begin{equation}\label{2varEqual1}
\Var Z_{\nu'_{\rho,\lambda}} (t) = \frac{2t}{\lambda \rho(\lambda - \rho)} + O(t^{2/3}) \,.   
\end{equation}

\end{demo}



\subsection{The invariant two-class Hammersley process}
Recall the construction of the invariant two-class Hammersley process $(M_{\nu_{\rho}}^{t}, \bar{M}^t_{\xi'})_{t\geq 0} $, described in Section \ref{twoclass}, and its connection with the two-line process and $M/M/1$ queues. Next, we use this construction to show Theorem \ref{theo:ZZ}.\\

\begin{demo}{\bf Proof of Theorem \ref{theo:ZZ}.}
Consider the Hammersley process with two classes of particles $(M_{\nu_{\rho}}^{t}, \bar{M}^t_{\xi'})_{t\geq 0} $ with invariant initial measure conditioned to have a second class particle $Z_{\xi'}$ at the origin, and constructed as a function of the multiline process $(M_{\alpha_1}^{t},M_{\alpha_2}^{t}, \bP)$. Let $(Q^{0}(x))_{x \in \RR}$ be the $M/M/1$ queue generated by $\alpha_1$ and $\alpha_2$ which gives origin to the initial measure $(M_{\nu_{\rho}}^{0},\bar{M}^t_{\xi'})$. As we conditioned to have a second class particle at the origin, $Q^{0}(0_{+}) = 0$. From the point of view of $Z_{\xi'}$, we can consider all the particles to its right as first class particles. By the construction of the queue, $Z_{\xi'}(0)$ sees a Poisson process of rate $\lambda$ to its right (process $\alpha_2$). To the left of $Z_{\xi'} (0)$, only first class particles affect its trajectory, that is, if we add or subtract second class particles located to the left of $Z_{\xi'} (0)$ its trajectory simply does not change. From the queue point of view, $Z_{\xi'}(0)$ only sees the effective services of the queue, which is empty at the origin. This is not a Poisson process but we will show that it is very similar to a Poisson process of density $\rho$. \\

We construct a new queue $(G^{0}(x))_{x \in \RR}$ whose arrival process and service processes are also $\alpha_1$ and $\alpha_2$, but $G^0 (0_{+})$ is a geometric random variable with parameter $\frac{\rho}{\lambda}$, independent of $\alpha_1$ and $\alpha_2$. This new queue is stationary and thus satisfies Burke's theorem. This means that the effective services of $G^0$ form a Poisson process of rate $\rho$. Now, suppose $G^0 (0_{+}) = 0$, then $G^0 (x) = Q^0 (x)$ for all $x \in \RR$ and the output process of $(G^{0}(x))_{x \leq 0}$ (with decreasing time) coincides with $M_{\nu_{\rho}}^{0} \cup \bar{M}^0_{\xi'} $. Let $k \in \NN$, if $G^0 (0_{+}) = k+1$, then for $x \leq 0$, $G^0 (x) > Q^0 (x)$ when $x > Z_{-k} (0)$ and $G^0 (x) = Q^0 (x)$ when $x \leq Z_{-k} (0)$, where $Z_{-k} (0)$ denote the position at time $0$ of the $k$-th second class particle to the left of the origin in the two-class process. Therefore, $M_{\nu_{\rho}}^{0} (- \infty, Z_{-k} (0))$ is equal to the Poisson process with rate $\rho$ generated from the effective services of $(G^0 (x))_{x \leq 0}$ for any $k \in \NN$. \\ 

The second class particles are renewal points for the invariant two-class process, which implies that the process as seen from the second class particle is in equilibrium \cite{FeMaH}. Thereby, for all $t$ we can construct the measure as seen from $Z_{\xi'} (t)$ from a queue $(Q^{t}(x))_{x \in \RR}$ generated by $M_{\alpha_1}^{t}$ and $M_{\alpha_2}^{t}$ and couple with the output of a stationary $M/M/1$, $(G^{t}(x))_{x \in \RR}$, in the same way as before, and we obtain that $Z_{\xi'} (t) - Z_{-k}(t)$ has the same law as $|Z_{-k}(0)|$. Thus, when $x$ goes to $-\infty$, the density to the left of $Z_{\xi'} (t)$ goes to the same density to the left of $Z_{-k}(0)$, which is $\rho$ (uniformly in $t$, by the invariance in distribution). Therefore $Z_{\xi'} (t)$ is a microscopic shock associated to $(M_{\nu_{\rho}}^{t}, \bar{M}^t_{\xi'})_{t\geq 0}$. \\ 

Now, we proceed with the coupling between $Z_{\nu'_{\rho,\lambda}} (t)$ and $Z_{\xi'} (t)$. We start constructing the measure $\nu_{\rho,\lambda}^{'}$ in the following way: as before, to the left of the origin it coincides with the used services of queue $(G^{0}(x))_{x < 0}$, and to the right of the origin it coincides with the service times of queue $(G^{0}(x))_{x > 0}$, that is, it coincides with $\alpha_2$. At the origin we put a second class particle $Z_{\nu'_{\rho,\lambda}} (t)$ independently of $G^0 (0_{+})$ (that is, if $G^0 (0_{+}) > 0$, we substitute the first class particle at the origin by a second class particle). We run the classical Hammersley process $(M_{\nu_{\rho,\lambda}^{'}}^t)_{t\geq 0}$ using the same jumping process $\bP$ of the two-class stationary process $(M_{\nu_{\rho}}^{t}, \bar{M}^t_{\xi'})_{t\geq 0} $. From this coupling we get that, for $x>0$, 
\begin{equation}\nonumber
 |\nu_{\rho,\lambda}^{'}(-x)|  \geq |M_{\nu_{\rho}}^{0} (-x)|\,\,\mbox{ and }\,\, |\nu_{\rho,\lambda}^{'}(x)|  = \alpha_2 (x)\,,
\end{equation}
which implies that,
\begin{equation}\nonumber
  Z_{\nu'_{\rho,\lambda}} (t) \leq Z_{\xi'} (t)\,,\mbox{ for all $t\geq 0$ }\,.
\end{equation}

Now, suppose that $G^0 (0_{+}) = k+1$ (the case $G^0 (0_{+}) = 0$ is straightforward). We add a new second class particle $Z_G (t)$ to the process $(M_{\nu_{\rho,\lambda}^{'}}^t)_{t\geq 0}$ with initial position $Z_G (0) =: Z_{-k-1} (0)$, where $Z_{-k-1} (0)$ is the initial position of the $(k+1)$-th second class particle to the left of the origin in the two-class process $(M_{\nu_{\rho}}^{t}, \bar{M}^t_{\xi'})_{t\geq 0} $. Note that $Z_G (t)$ does not affect the trajectory of $Z_{\nu'_{\rho,\lambda}} (t) $ like any second class particle to its left. Also, $Z_G (t) = Z_{-k-1} (t)$ for all $t\geq 0$. To see this, remind we have $M_{\nu_{\rho,\lambda}^{'}}^0 (- \infty, Z_{-k} (0)) = M_{\nu_{\rho}}^{0} (- \infty, Z_{-k} (0))$ and $M_{\nu_{\rho,\lambda}^{'}}^0 (Z_{-k-1} (0), +\infty) = M_{\nu_{\rho}}^{0} \cup \bar{M}^t_{\xi'} (Z_{-k-1} (0), +\infty)$, which shows the initial condition of the process $(M_{\nu_{\rho,\lambda}^{'}}^t)_{t\geq 0}$ as seen from $Z_G (t)$ is equal to the initial condition of the process $(M_{\nu_{\rho}}^{t}, \bar{M}^t_{\xi'})_{t\geq 0} $ as seen from $Z_{-k-1} (t)$. By the coupling, both use the same jump process $\bP$, which shows $Z_G (t) = Z_{-k-1} (t)$ for all $t\geq 0$. Since  $Z_{\nu'_{\rho,\lambda}} (t)  \geq Z_G (t)$ it follows that 
$$ Z_{-k-1} (t) \leq Z_{\nu'_{\rho,\lambda}} (t)  \leq Z_{\xi'} (t)\,,$$
from where we conclude that 
$$|Z_{\nu'_{\rho,\lambda}} (t)  - Z_{\xi'} (t) | \leq J_t:=  Z_{\xi'} (t) - Z_{-k-1} (t)\,{\buildrel \iD \over =}\,|Z_{-k-1} (0)|\,$$
(since the two-class process is in equilibrium). Therefore, $Z_{\nu'_{\rho,\lambda}} (t) $ is a microscopic shock associated to the process $(M_{\nu_{\rho,\lambda}^{'}}^{t})_{t\geq 0}$, and the random variables $(J_t)_{t\geq 0}$ are identically distributed with the same law of $|Z_{-k-1} (0)|$, which has finite moments. 

\end{demo}


\vspace{1.0cm}


\begin{thebibliography}{10}

\bibitem{AD} Aldous, D.\ and\ Diaconis, P.: Hammersley's interacting particle process and longest increasing subsequences. Probab. Theory Related Fields {\bf 103}, no.~2, 199--213 (1995). 



\bibitem{BaRa} Baik, J.\ and\ Rains, E. M.: Limiting distributions for a polynuclear growth model with external sources. J. Statist. Phys. {\bf 100}, no.~3-4, 523--541 (2000). 



\bibitem{CG1} Cator, E.\ and\ Groeneboom, P.: Hammersley's process with sources and sinks. Ann. Probab. {\bf 33}, no.~3, 879--903 (2005). 


\bibitem{CG2}  Cator, E.\ and\ Groeneboom, P.: Second class particles and cube root asymptotics for Hammersley's process. Ann. Probab. {\bf 34}, no.~4, 1273--1295 (2006). 


\bibitem{CPS} Cator, E. A.,  Pimentel, L. P. R.\ and\ Souza, M. W. A.: Influence of the initial condition in equilibrium last-passage percolation models. Electron. Commun. Probab. {\bf 17}, no. 7, 7 pp (2012). 


\bibitem{Fe} Ferrari, P. A.: Shocks in the Burgers equation and the asymmetric simple exclusion process. in {\it Statistical physics, automata networks and dynamical systems (Santiago)}, 25--64, Math. Appl., 75 Kluwer Acad. Publ., Dordrecht (1990). 


\bibitem{Fe2} Ferrari, P. A.: Shock fluctuations in asymmetric simple exclusion. Probab. Theory Related Fields {\bf 91}, no.~1, 81--101 (1992). 


\bibitem{FeFoS}  Ferrari, P. A.\ and\ Fontes, L. R. G.: Shock fluctuations in the asymmetric simple exclusion process. Probab. Theory Related Fields {\bf 99}, no.~2, 305--319 (1994).   



\bibitem{FeMaH}  Ferrari, P. A.\ and\  Martin,J. B.: Multiclass Hammersley-Aldous-Diaconis process and multiclass-customer queues.
\newblock \emph{Ann. Inst. Henri Poincar\'e Probab. Stat.} {\bf 45} , no.~1, 250--265 (2009). 


\bibitem{Ha}  Hammersley, J. M.: A few seedlings of research, in {\it Proceedings of the Sixth Berkeley Symposium on Mathematical Statistics and Probability (Univ. California, Berkeley, Calif., 1970/1971), Vol. I: Theory of statistics}, 345--394, Univ. California Press, Berkeley, CA (1972).


\bibitem{SeM}  Sepp\"al\"ainen, T.: A microscopic model for the Burgers equation and longest increasing subsequences. Electron. J. Probab. {\bf 1}, no.\ 5, approx.\ 51 pp.\ (electronic) (1996).


\bibitem{SeD}  Sepp\"al\"ainen, T.: Diffusive fluctuations for one-dimensional totally asymmetric interacting random dynamics. Comm. Math. Phys. {\bf 229}, no.~1, 141--182 (2002). 


\end{thebibliography}
\end{document}